\documentclass[12pt,twoside]{article} 
\usepackage{amssymb,amsfonts,mathrsfs,graphicx}
\usepackage[small,sc]{caption} 

\newtheorem{theorem}{Theorem}[section]

\newtheorem{proposition}[theorem]{Proposition}
\newtheorem{corollary}[theorem]{Corollary}

\newtheorem{rmrk}[theorem]{Remark}

\makeatletter
\@addtoreset{equation}{section}
\makeatother

\newcommand{\figmod}[3] {
\medskip\smallskip
\begin{figure}[htb]
  \centering
  \includegraphics[width=#2]{img_#1.pdf}
  \begin{minipage}[t]{0.85\linewidth} 
    \vspace{-8pt}
    \caption{#3}
    \protect\label{#1}
  \end{minipage}
\end{figure}
\medskip
}

\newcommand{\figmpng}[3] {
\medskip\smallskip
\begin{figure}[htb]
  \centering
  \includegraphics[width=#2]{img_#1.png}
  \begin{minipage}[t]{0.85\linewidth} 
    \caption{#3}
    \protect\label{#1}
  \end{minipage}
\end{figure}
\medskip
}

\newenvironment{remark}
{\begin{rmrk} \em}
{\end{rmrk}}

\newcommand{\R} {\mathbb{R}}

\newcommand{\Z} {\mathbb{Z}}
\newcommand{\N} {\mathbb{N}}

\newcommand{\skippar} {\par\medskip}

\newcommand{\article}[3] {\textsc{{#1}}, {\itshape {#2}}, {{#3}}.}
\newcommand{\book}[3] {\textsc{{#1}}, {\itshape {#2}}, {{#3}}.}
\newcommand{\vol} {\textbf}
\newcommand{\eps} {\varepsilon}
\newcommand{\rset}[2] {\left\{ #1 \: \left| \: #2 \right. \! \right\} }

\newcommand{\om} {\omega}
\renewcommand{\l} {\ell}
\newcommand{\mM} {\mathcal{M}}
\newcommand{\mMa} {\mathcal{M}_\alpha}
\newcommand{\mN} {\mathcal{N}}
\newcommand{\tu} {\mathcal{T}}  
\newcommand{\ta} {Q}            
\newcommand{\si} {\mathcal{S}}  
\newcommand{\bo} {B}            
\newcommand{\ob} {O}            
\newcommand{\mD} {\mathcal{D}}
\newcommand{\tnl} {T_{\mN : \l}}
\newcommand{\sia} {\si_\alpha}

\newcommand{\bi} {billiard}
\newcommand{\fn} {function}
\newcommand{\me} {measure}
\newcommand{\tr} {trajector}
\newcommand{\erg} {ergodic}
\newcommand{\sy} {system}
\newcommand{\hyp} {hyperbolic}
\newcommand{\pr} {probability}
\newcommand{\ra} {random}
\newcommand{\dsy} {dynamical system}
\renewcommand{\o} {orbit}

\begin{document}

\title{\textbf{Recurrence and higher ergodic properties \\ 
for quenched random Lorentz tubes in \\
dimension bigger than two}}

\author{\scshape 
Marcello Seri\,%
\thanks{
Dipartimento di Matematica, Universit\`a di Bologna, 
Piazza di Porta San Donato 5, 40126 Bologna, Italy. E-mail:
\texttt{seri,lenci,desposti,cristadoro@dm.unibo.it}.
}
\thanks{
Department of Mathematics, University of Erlangen-Nuremberg,
Bismarckstr.\ 1 1/2, 91053 Erlangen, Germany.}
,
Marco Lenci\,$^*$ \\
\scshape 
Mirko degli Esposti\,$^*$, Giampaolo Cristadoro\,$^*$
}

\date{Final version for J.~Stat.~Phys., June 2011}

\maketitle

\begin{abstract}
  We consider the billiard dynamics in a non-compact set of $\R^d$
  that is constructed as a bi-infinite chain of translated copies of
  the same $d$-dimensional polytope. A random configuration of
  semi-dispersing scatterers is placed in each copy. The ensemble of
  dynamical systems thus defined, one for each global realization of
  the scatterers, is called \emph{quenched random Lorentz tube}. Under
  some fairly general conditions, we prove that every system in the
  ensemble is hyperbolic and almost every system is recurrent,
  ergodic, and enjoys some higher chaotic properties.

\end{abstract}

\section{Introduction}
\label{sec-intro}

A $d$-dimensional Lorentz tube (LT) is a Lorentz gas, in Euclidean
$d$-space, that is confined to a subset $\tu$ which is infinitely
extended in one dimension.

As a prototype, think of an infinite square-section cylinder in
$\R^3$, in whose interior a countable number of convex scatterers are
placed approximately with the same density (see
Fig.~\ref{hdtubes_hole}). A material point travels inertially in the
free region of $\tu$, until it collides with either a scatterer or the
boundary of $\tu$ (from now on, the latter will be referred to as a
scatterer as well). Assuming the scatterer to be infinitely massive,
the collision is totally elastic, i.e., the outgoing velocity $v^+$ is
derived from the incoming velocity $v^-$ by reversing the normal
component of $v^-$ relative to the plane of collision.

\figmod{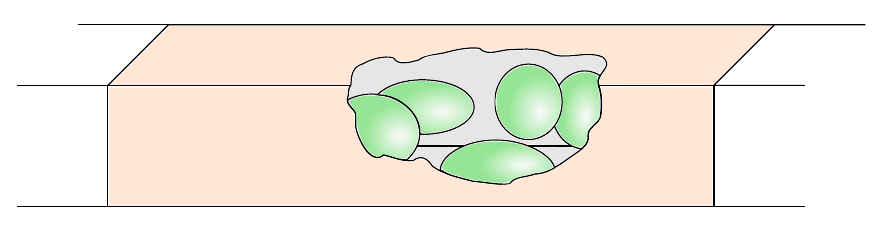} {13cm} {An example of a 3D LT.}

In the terminology of \dsy s, a \sy\ like this is an \emph{extended
  semi-dispersing \bi}. The term `extended' refers to the fact that
the configuration space is not compact and the relevant physical \me\
on it is infinite. Also, it is a semi-dispersing \bi\ because the
particle undergoes a \bi-like dynamics with bouncing walls that are
either flat or convex, as seen from the particle. (We invite the
reader to avoid confusion with \bi s that are semi-dispersing because
the sectional curvature of the scatterers can be either positive or
zero; e.g., 3D billiards with cylindrical scatterers.)

It is a celebrated fact that semi-dispersing \bi s give rise to
chaotic dynamics \cite{cm}---using the term `chaotic' in a very lax
sense here---and so it seems sound to use models like this to study
the motion of small particles (e.g., electrons) in thin wires, in
whose interior a configuration of obstacles (e.g., atomic nuclei)
makes the motion chaotic. (This was more or less Lorentz's original
motivation \cite{lo}; cf.\ also \cite{kf, aacg, lwwz, heal} and
references therein.)

\skippar

This note is a follow-up to an article that three of the present
authors have published recently \cite{cls}, where these types of \sy s
are studied in two dimensions. We refer the reader to the introduction
of that paper for a better description of the physical and
mathematical motivations behind this research.

Here we just outline our main result and its consequences in terms of
the stochastic properties of the dynamics we consider. Our chief
interest, as far as this note is concerned, is in the recurrence of
these types of \sy s. Recurrence is the most basic property one needs
to establish of extended \sy s in order to study their chaotic
properties.  (It is hard to claim that a certain dynamics
``randomizes'' the state of the \sy---more precisely, decorrelates it
from its initial condition---if a non-negliglible part of the phase
space is made up of \tr ies that escape to infinity, thus giving no
asymptotic contribution to the state of the \sy\ in any given compact
region.)

In fact, for our LTs, we will see below that recurrence is a
sufficient condition for a number of stronger \erg\ properties.

For these effectively one-dimensional \sy s one would expect
recurrence to be a typical property. To make this point, we define a
fairly large and representative \me d family (in the language of
statistical mechanics, an \emph{ensemble}) of LTs and ask if the
typical element is recurrent, in the sense of Poincar\'e (which
coincides with the intuitive meaning of the word here).

The family is defined roughly as follows: The tube $\tu$ is made up of
a countable number of congruent $d$-dimensional polytopes (henceforth
\emph{cells}) having two parallel and congruent facets (henceforth
\emph{gates}), whereby each cell is attached to its two adjacent cells
(Fig.~\ref{hdtubes_hole} is an example of this). In each cell we put a
\ra\ configuration of convex scatterers, according to a very general
\pr\ law.  Each global realization of scatterers defines a different
LT. This type of structure, in which one \ra ly chooses a \dsy\ and
then follows its deterministic dynamics, is called `quenched \ra\
\dsy'.  So what we have is a \emph{quenched \ra\ LT}.

If the \sy\ verifies some geometric conditions (most of which are
rather general, some less), we prove that almost surely, in the sense
of the \pr, the LT is recurrent. Since it can be proved that a
recurrent LT is also \erg, and that suitable first-return maps are
strongly chaotic (at least $K$-mixing \cite{aa}), an important
corollary of our work is that the typical LT in our family is chaotic.

As already mentioned, the main contribution of this paper is an
extension of the results of \cite{cls} to dimension $d > 2$. It is
known by experts in the field that semi-dispersing \bi s in dimension
three and higher present specific subtleties and difficulties.
Therefore, we made an effort to detail the parts of the proofs that
deal with such difficulties, while giving a looser exposition of the
remaning arguments (as they can be found elsewhere as well). The paper
is organized as follows: Section \ref{sec-pre}, which should be
accessible to the reader with a minimal background in \dsy s, contains
the precise formulation of the results. Section \ref{sec-sketch} gives
an outline of the main proof, partly referring to previous work by
some of the present authors. In Section \ref{sec-detail}, which is the
most technical part of the paper, we give precise proofs for the
arguments that are specific to the \sy s at hand.

\bigskip
\noindent
\textbf{Acknowledgments.}\ We thank Gianluigi Del Magno and Domokos
Sz\'asz for useful discussions. This work was partially supported by
the FIRB-``Futuro in Ricerca'' Project RBFR08UH60 (MIUR, Italy).

\section{Mathematical formulation of the results}
\label{sec-pre}

Consider a closed $d$-polytope $C_0 \in \R^d$ that has two parallel
and congruent facets. Denoting said facets $G^1$ and $G^2$, call
$\tau$ the translation of $\R^d$ that takes $G^1$ into $G^2$, and
define, for $n \in \Z$, $C_n := \tau^n (C_0)$. Each $C_n$ is called a
\emph{cell} and $\tu := \bigcup_{n \in \Z} C_n$ is called the
\emph{tube}, see Fig.~\ref{hdgates}.

\figmod{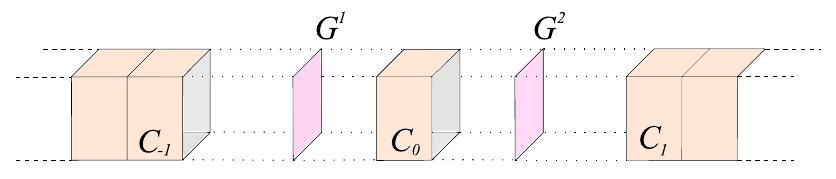} {13cm} {Assembling an LT as a chain of cells with
  gates.}

For each $n$, a family of closed, pairwise disjoint, piecewise smooth,
convex sets $\ob_{n,i} \subset C_n$ ($i = 1, \ldots, N$) is given. We
refer to this family as the \emph{local} configuration of scatterers
in the cell $C_n$ (note that some $\ob_{n,i}$ might be empty, so
different cells might have a different number of scatterers). This
configuration is \ra, in the sense that each $\ob_{n,i} =
\ob_{n,i}(\l_n)$ is a \fn\ of the random parameter $\l_n \in \Omega$,
where $\Omega$ is some \me\ space. The sequence $\l := ( \l_n )_{n \in
  \Z} \in \Omega^\Z$, which thus describes the \emph{global}
configuration of scatterers in the tube $\tu$, is a stochastic process
obeying the \pr\ law $\Pi$, whose properties are given later, cf.\ in
particular (A1).

For each realization $\l$ of the process, we consider the \bi\ in the
\emph{table} $\ta_\l := \tu \setminus \bigcup_{n \in \Z}
\bigcup_{i=1}^N \ob_{n,i}(\l_n)$. This is the \dsy\ $(\ta_\l \times
S^{d-1}, \phi_\l^t, m_\l)$, where $S^{d-1}$ is the unit sphere in
$\R^d$ and $\phi_\l^t : \ta_\l \times S^{d-1} \longrightarrow \ta_\l
\times S^{d-1}$ is the \emph{\bi\ flow}, whereby $(q_t, v_t) =
\phi_\l^t (q, v)$ represents the position and velocity at time $t$ of
a point particle with initial conditions $(q, v)$, undergoing free
motion in the interior of $\ta_\l$ and Fresnel collisions at $\partial
\ta_\l$, i.e., if $q_t \in \partial \ta_\l$, then
\begin{equation}
  \label{coll}
  v_{t^+} =  v_{t^-} - 2 ( v_{t^-} \cdot o_{q_t} ) \, o_{q_t},
\end{equation}
where $o_{q_t}$ is the inner unit normal to $\partial \ta_\l$ at
$q_t$.  (Notice that in this Hamiltonian \sy\ the conservation of
energy corresponds to the conservation of speed, which is thus
conventionally fixed to 1.) Lastly, $m_\l$ is the Liouville invariant
\me\ which, as is well known, is the product of the Lebesgue \me\ on
$\ta_\l$ and the Haar \me\ on $S^{d-1}$.

We call this \sy\ the \emph{LT corresponding to the realization $\l$},
or simply the \emph{LT $\l$}. As $\l$ ranges in the \pr\ space
$(\Omega^\Z, \Pi)$, we have a \ra\ family, or an ensemble, of \dsy
s. This structure is referred to as a `quenched \ra\ \dsy'. As we
shall see later, the situation is simplified by the fact that these
\dsy s can be reformulated in such a way that they all share the same
phase space and invariant \me.

In the remainder, whenever there is no risk of ambiguity, we drop the
dependence on $ \l$ from all the notation. Also, we call
\emph{universal constant} any bound that depends on none of the
quantities explictly or implicitly involved in the inequality at hand,
in particular on $\l$.

\skippar

We assume the following:

\begin{itemize}
\item[(A1)] $\Pi$ is \erg\ for the left shift $\sigma: \Omega^\Z
  \longrightarrow \Omega^\Z$.

\item[(A2)] There exists a universal constant $K_1 \in \Z^+$ such that
  (for all realizations $\l \in \Omega^\Z$) $\partial \ob_{n,i}$ is
  made up of at most $K_1$ compact, connected, uniformly $C^3$
  (w.r.t.\ $n,i$) subsets of algebraic varieties (SSAVs), which may
  intersect only at their borders. These borders, which thus have
  codimension larger than one, will be generically referred to as
  \emph{edges}.

\item[(A3)] If $q$ is a smooth point of $\partial \ta$, let
  $\mathbf{k}(q)$ be the second fundamental form of $\partial \ta$ at
  $q$. There are two universal constants $k_M > k_m > 0$ such that,
  for all smooth $q \in \partial \ta$, either the SSAV that $q$
  belongs to is a piece of a hyperplane or
  \begin{displaymath}
    k_m \le \mathbf{k}(q) \le k_M,
  \end{displaymath}
  where the inequalities are meant in the sense of the quadratic
  forms.

\item[(A4)] There exist universal constants $L > 0$, $K_3 \in \Z^+$,
  and $\eta \in (0, \pi/2)$ such that, in each portion of \tr y of
  length (equivalently, duration) $L$, there are at most $K_3$
  collisions; and at least one collision with a \emph{dispersing}
  (i.e., non-flat) part of $\partial \ta$ and such that the angle of
  incidence (relative to the normal at the collision point) is less
  than $\pi/2 - \eta$. Notice that the above implies the so-called
  \emph{finite-horizon condition}: the free flight is bounded above.

\item[(A5)] A singular \tr y is a \tr y which has tangential
  collisions or collisions with the edges of $\partial \ta$ (in which
  case it conventionally ends there). Using this terminology, we
  assume that, for a.e.\ $\l$ and all $i,j \in \{ 1,2 \}$, there is a
  non-singular \tr y entering $C_0$ through $G^i$ and leaving it
  through $G^j$.
\end{itemize}

The next and last assumption has to do with the well-known fact that a
semi-dispersing \bi\ is a discontinuous (and indeed singular) \dsy. It
will be formulated in full mathematical rigor in Section
\ref{sec-detail}, after the necessary definitions are given. Here we
give a descriptive version which will be quickly understood by the
``\hyp\ \bi ist''.  We anticipate, however, that this assumption is
verified for a reasonable class of perturbations of a periodic LT or
in the case in which $\Omega$ is finite, that is, the local
configuration in each cell is chosen from a finite number of
possibilities.

\begin{itemize}
\item[(A6)] There exist a universal constant $K_4 > 0$ such that the
  Lebesgue \me\ of the $\delta$-neighborhood of each smooth piece of
  the singularity set does not exceed $K_4 \, \delta$.
\end{itemize}

We then have:

\begin{theorem}
  \label{thm-main}
  Under assumptions \emph{(A1)-(A6)}, the quenched \ra\ LT is almost
  surely recurrent, that is, for $\Pi$-a.e.\ $\l \in \Omega^\Z$, the
  LT $\l$ is Poincar\'e recurrent.
\end{theorem}

For the sake of completeness, we recall what Poincar\'e recurrence
means in our context: Given a measurable $A \subset \ta \times
S^{d-1}$, for $m$-a.e.\ $(q,v) \in A$, there is an unbounded sequence
of times $t_j$ such that $\phi^{t_j} \in A$.

\begin{remark}
  Although the theorem is valid in any dimension $d \ge 2$, the
  assumptions given earlier were designed for the case $d \ge 3$: in
  the two-dimensional case, as presented in \cite{cls}, the hypotheses
  are substantially weaker.
\end{remark}

Theorem \ref{thm-main} has deep implications. For a fixed $\l$, let
$D$ be a finite union of the SSAVs that make up $\partial \ta$; cf.\
(A2). (So, for example, $D$ could be one smooth piece of some
$\partial \ob_{n,i}$, or the whole of it.) Denote
\begin{equation}
  \mD := \rset{(q,v) \in D \times S^{d-1}} {v \cdot o_q \ge 0},
\end{equation}
that is, $\mD$ is the submanifold in phase space corresponding to the
post-collisional position-velocity pairs (henceforth \emph{line
  elements}) based in $D$. It is apparent that, if the LT $\l$ is
recurrent, the first-return map onto $\mD$ is well-defined almost
everywhere w.r.t.\ the natural \me\ on $\mD$ (see below): we call this
map $T_\mD$.  The Liouville \me\ $m$ induces a $T_\mD$-invariant \me\
on $\mD$, which we denote $\mu_\mD$. It is well known that $d\mu_\mD
(q,v) = (v \cdot o_q) dq dv$, where $dq$ is the volume element in
$\partial \ta$ and $dv$ is the volume, or Haar, element in $S^{d-1}$
\cite{cm}. It is a consequence of (A2) that $D$ has a finite volume
so, upon normalization, we may assume that $\mu_\mD (\mD) = 1$.

To avoid misunderstandings, let us recall that both \dsy s $(\ta
\times S^{d-1}, \phi^t, m)$ and $(\mD, T_\mD, \mu_\mD)$ depend on the
choice of $\l$---the subscript has been removed only to lighten the
notation.

An important result is the following:

\begin{theorem}
  \label{thm-erg}
  If the LT $\l$ is recurrent, then $(\ta \times S^{d-1}, \phi^t, m)$
  is \erg. Also, for any choice of $D$ as described above, $(\mD,
  T_\mD, \mu_\mD)$ is K-mixing (thus mixing and \erg).
\end{theorem}

Combining Theorems \ref{thm-main} and \ref{thm-erg}, we obtain:

\begin{corollary}
  \label{cor-erg}
  $(\ta \times S^{d-1}, \phi^t, m)$ is \erg\ and $(\mD, T_\mD,
  \mu_\mD)$ is K-mixing for $\Pi$-a.e.\ choice of $\l \in \Omega^\Z$.
\end{corollary}

Since generating an example of an LT that verifies (A1)-(A6) may not
be immediate, we present one in Appendix \ref{app}.

\section{Flow of the proofs}
\label{sec-sketch}

The proofs of Theorems \ref{thm-main} and \ref{thm-erg} follows
exactly the same strategy as the corresponding proofs in \cite{cls,
  l2, l1}. For the convenience of the reader, though, we are going to
outline them in this section, with particular regard to Theorem
\ref{thm-main}.  Some of the intermediate results will present
complications due to the higher-dimensional setting. We will explain
how to prove these results in closer detail in
Section~\ref{sec-detail}.

The reader is warned that there are a few changes of notation compared
to \cite{l1, l2, cls}.

The first step of the proof of Theorem \ref{thm-main} consists in
showing the \hyp\ properties of \emph{each} \dsy\ in the ensemble. So,
for the time being, we fix $\l \in \Omega^\Z$ and describe the LT $\l$
by means of a certain Poincar\'e map which we introduce
momentarily. For $n \in \Z$, denote by $\{ \bo_{n,j}
\}_{j=1}^{\kappa_n}$ the collection of all the \emph{dispersing}
pieces of boundary in the cell $C_n$ (it is important that $\bo_{n,j}$
be a whole SSAV among those mentioned in (A2), and not just a portion
of it).  From our hypotheses, $\kappa_n$ is bounded above by the
universal constant $K_1 N$.

Henceforth, we will indicate the index $(n,j)$ with the symbol
$\alpha$, and the space of all such indices with $\mathcal{A}$.
Evidently, $\mathcal{A}$ is countable.  For $\alpha \in \mathcal{A}$,
define
\begin{equation}
  \label{def-malpha}
  \mMa : = \rset{(q,v) \in \bo_\alpha \times S^{d-1}} {v \cdot
  o_q \ge \eps},
\end{equation}
where $\eps := \cos(\pi/2 - \eta)$ and $\eta$ is the universal
constant that appears in (A4). (Again, $o_q$ is the inner unit normal
to $\partial \ta$ at the point $q$.)

All the above definitions clearly depend on $\l$: let us now reinstate
this dependance in the notation and denote $\mM_\l := \bigcup_{\alpha
  \in \mathcal{A}} \mMa$. By (A4), $\mM_\l$ is a global cross-section
for the flow $\phi_\l^t$. We call $T_\l$ and $\mu_\l$, respectively,
its Poincar\'e map and the invariant \me\ induced by the Liouville
\me\ $m_\l$ on $\mM_\l$ (of course, up to a constant factor, $\mu_\l$
has the same density as the \me\ $\mu_\mD$ introduced in Section
\ref{sec-pre}).  Each \dsy\ thus defined possesses some basic \hyp\
and \erg\ properties which we now outline, in a rather undetailed way.

\begin{theorem}
  \label{thm-ltl}
  The following holds for the \dsy\ $(\mM_\l, T_\l, \mu_\l)$:
  \begin{itemize}
  \item[(a)] The \sy\ is uniformly \hyp\ w.r.t.\ the natural metric in
    $\mM_\l$.

  \item[(b)] There is a \hyp\ structure, in the sense that there
    exists local stable and unstable manifolds (LSUMs) at a.e.\ point
    of $\mM_\l$. Also, the two corresponding (invariant) foliations,
    when \me d with a Lebesgue-equivalent $(d-1)$-dimensional \me, are
    absolutely continuous relative to $\mu_\l$.
  
  \item[(c)] Local \erg ity holds. This means that, for each $\alpha$,
    a.a.\ pairs of points in $\mMa$ are connected by a chain of
    alternating LSUMs that intersect transversally. The intersection
    points can be chosen out of a predetermined full-\me\ subset of
    $\mMa$.
  \end{itemize}
\end{theorem}

In the second step of the proof we represent the LT $\l$ in yet
another way, which will be more convenient later on. What we do is, we
introduce a different cross section for the same flow. For $n \in \Z$
and $j \in \{1,2\}$, denote by $G_n^j := \tau^n (G^j)$ the two
\emph{gates} to the cell $C_n$ and by $o_j$ be the inner normal to
$G_n^j$, relative to $C_n$ (thus $o_2 = -o_1$). Consider
\begin{equation}
  \label{def-nnj}
  \mN_n^j : = \rset{(q,v) \in G_n^j \times S^{d-1}} {v \cdot
  o_j > 0},
\end{equation}
that is, the collection of (almost) all line elements entering $C_n$
from the ``left'' or from the ``right'', depending on $j$. The global
cross section that we use this time is $\mN := \bigcup_{n \in \Z} \,
\bigcup_{j=1,2} \mN_n^j,$ while the corresponding Poincar\'e map we
denote $\tnl$.

The gates are sometimes referred to as \emph{transparent walls},
because, in the theory of \bi s, the corresponding map has virtually
the same properties as an ordinary \bi\ map, such as $T_\mD$ or
$T_\l$. In particular, it preserves a \me\ $\mu_\mN$ that has the same
\fn al form as the \me s $\mu_\mD$ and $\mu_\l$.

We end up with the triple $(\mN, \tnl, \mu_\mN)$. Notice that neither
the phase space nor the \me\ depend on $\l$, which is precisely what
makes this \dsy\ more convenient than the one previously introduced
(and gives further justification as to why the whole ensemble is
called `quenched \ra\ \dsy': we have a family of maps that are defined
on the same space and preserve the same \me).

At this point, one might ask why the map $T_\l$ was defined at
all. The reason is, we needed to prove Theorem \ref{thm-ltl} first, in
order to obtain the corresponding results for $(\mN, \tnl, \mu_\mN)$.
In fact, it is not hard to verify that the latter \sy\ inherits the
\hyp\ structure of the former: One constructs the local stable
manifolds (LSMs) and local unstable manifolds (LUMs) of $(\mN, \tnl,
\mu_\mN)$ as push-forwards, respectively pull-backs, of the LSMs and
LUMs of $(\mM_\l, T_\l, \mu_\l)$ \cite{l2}. This is possible because
the first \sy\ has fewer singular \tr ies than the second so, for
example, when a LSM of the second \sy\ is pushed forward by the flow,
no cuts occur at all positive times---ensuring that the defined
push-forward is indeed a LSM.

It is also rather easy to check that uniform \hyp ity is maintained.
The result that we are mostly interested in, however, is the analog of
Theorem \ref{thm-ltl}\emph{(c)}:

\begin{theorem}
  \label{thm-ltl2}
  The \dsy\ $(\mN, \tnl, \mu_\mN)$ is locally \erg\ in the following
  sense: For any $n \in \Z$ and $j \in \{ 1,2 \}$, a.a.\ pairs of
  points in $\mN_n^j$ are connected by a chain of alternating LSUMs
  that intersect transversally. The intersection points can be chosen
  out of a predetermined full-\me\ subset of $\mN_n^j$.
\end{theorem}

In the third step of the proof we use the so-called `point of view of
the particle' (PVP). It consists of a \emph{finite-\me} \dsy\ that,
together with a suitable observable, describes the dynamics of
\emph{all} the \o s in \emph{all} the realizations of the LT.  The
idea is that, instead of following a given \o\ from one cell to
another, with every iteration of the dynamics we shift the LT in the
direction opposite to the \o's displacement, so that the point always
lands in the same cell (conventionally $C_0$).  We briefly outline the
construction of this \dsy, referring the reader to \cite{cls, l2} for
more detailed explanations.

Let $\mN_0 := \mN_0^1 \cup \mN_0^2$ be the cross-section corresponding
to the gates of $C_0$, and $\mu_0$ the normalized \bi\ \me\ on it.
For a given $\om \in \Omega$, determining the local configuration in
$C_0$, define a map $R_\om: \mN_0 \longrightarrow \mN_0$ as follows.
Trace the forward \tr y of $(q,v) \in \mN_0$ until it crosses $G^1$ or
$G^2$ for the first time (almost all \tr ies do). This occurs at a
point $q_1$ with velocity $v_1$. If, for $\epsilon \in \{-1,+1\}$,
$C_{\epsilon}$ is the cell that the particle enters upon leaving
$C_0$, define
\begin{equation}
  \label{eq-def-R}
  R_\om (q,v) := (\tau^{-\epsilon} (q_1), v_1).
\end{equation}
Clearly $R_\om(q,v) \in \mN_0$, and $R_\om$ preserves $\mu_0$ for
every $\om$. Next define the so-called \emph{exit function} $e: \mN_0
\times \Omega \longrightarrow \{-1,+1\}$ via the formula $e(q,v; \,
\om) := \epsilon$. From now on, we indicate line elements $(q,v)$ with
the letter $x$.

The PVP \sy\ $(\Sigma, F, \lambda)$ is defined by: 

\begin{itemize}
\item $\Sigma := \mN_0 \times \Omega^\Z$.

\item $F(x,\l) := (R_{\l_0} (x), \sigma^{e(x,\l_0)} (\l))$, which
  defines a map $\Sigma \longrightarrow \Sigma$. Here $\l_0$ is the
  $0^\mathrm{th}$ component of $\l$ and $\sigma$ is the left shift on
  $\Omega^\Z$, introduced in (A1) (therefore $\sigma^\epsilon (\l) =
  \{\l'_n\}_{n\in\Z}$, with $\l'_n := \l_{n + \epsilon}$).

\item $\lambda := \mu_0 \times \Pi$. Clearly, $\lambda (\Sigma) =
  1$. Also, using that $F$ is invertible, $R_\om$ preserves $\mu_0$,
  and $\sigma$ preserves $\Pi$, it can be seen that $F$ preserves
  $\lambda$.
\end{itemize}

Now, tolerating the abuse of notation whereby $e(x,\l) = e(x, \l_0)$,
let us think of the exit \fn\ $e$ as an integer-valued observable of
the \dsy\ just defined. We are interested in its \emph{cocycle}
(namely, Birkhoff sum) $\{ S_n \}_{n \in \N}$, given by $S_0(x,\l)
\equiv 0$ and
\begin{equation}
  S_n(x,\l) := \sum_{k=0}^{n-1} (e \circ F^k) (x,\l).
\end{equation}
A discrete cocycle, such as $\{ S_n \}$, is said to be
\emph{recurrent} if, for a.e.\ $(x,\l)$, there exists a subsequence
$\{ n_j \}$ such that $S_{n_j}(x,\l) = 0$, for all $j$. For
one-dimensional (i.e., $\Z$-valued) cocycles, a sufficient condition
for recurrence has long been known (see, e.g., \cite{at}):

\begin{proposition}
  \label{prop-rec}
  If $(\Sigma, F, \lambda)$ is \erg, and $e: \Sigma \longrightarrow
  \Z$ is integrable with $\int_\Sigma e = 0$, then the corresponding
  cocycle is recurrent.
\end{proposition}

(A beautiful $d$-dimensional version of this result was given by
Schmidt \cite{s}---see also a generalization in the Appendix of
\cite{l3}.)  It is not hard to verify that the recurrence of $\{ S_n
\}$ implies Theorem \ref{thm-main}. In fact, let us call a global
configuration $\l$ \emph{typical} if, for all $k \in \Z$ and
$\mu_0$-a.a.\ $x \in \mN_0$, $\{ S_n (x,\sigma^k(\l)) \}_{n \in \N}$
is recurrent. By the recurrence of the cocycle, Fubini's Theorem and
the denumerability of $\Z$, $\Pi$-a.e.\ $\l \in \Omega^\Z$ is
typical. On the other hand, by the above definition, in a typical LT
$\l$ almost all \o s come back to the cell where they started. Using
the Poincar\'e Recurrence Theorem on suitable first-return maps, one
easily checks that this is equivalent to the Poincar\'e recurrence of
$(\mN, \tnl, \mu_\mN)$ which, clearly, is the same as the recurrence
of $(\ta_\l \times S^{d-1}, \phi_\l^t, m_\l)$. In turn, since $\mM_\l$
is a global cross-section, that is equivalent to the recurrence of
$(\mM_\l, T_\l, \mu_\l)$ and implies the recurrence of $(\mD, T_\mD,
\mu_\mD)$.

Since $e$ is bounded and has zero average (this is clear by
time-reversal symmetry), in order to derive Theorem \ref{thm-main}
from Proposition \ref{prop-rec}, what remains to be shown is:

\begin{theorem}
  \label{thm-pvp}
  $(\Sigma, F, \lambda)$ is \erg.
\end{theorem}

This is proved precisely as Thm.~4.1 of \cite{cls}. The idea is to use
Theorem \ref{thm-ltl2} to show that each \emph{fiber} $\mN_0 \times
\{\l\}$ of $\Sigma$ is fully contained in one \erg\ component of the
\sy. In other words, the \erg\ decomposiontion is coarser than the
decomposition into fibers. Assumptions (A1) and (A5) then ensure that
a $\Pi$-typical fiber belongs to the \erg\ component of any other
typical fiber.

\skippar

We finish this section by giving a concise outline of the proof of
Theorem \ref{thm-erg}, following \cite{l1, l2}.

If $(\mN, \tnl, \mu_\mN)$ is recurrent, then the first-return map to
any $\mN_n^j$ is well-defined almost everywhere.  Theorem
\ref{thm-ltl2} ensures that it is also \erg, since a.a.\ points on a
LSUM belong to the same \erg\ component. This implies that no
$\mN_n^j$ can be split into two invariant sets of $\tnl$. In other
words, the \erg\ decomposition of $\tnl$ is coarser than the partition
of $\mN$ into connected components. But (A5) entails that $\tnl$
carries a positive \me\ of points from $\mN_n^j$ into $\mN_{n \pm
  1}^j$ and $\mN_n^{j \pm 1}$ (using the plus sign for $j=1$ and the
minus sign for $j=2$). This ensures that there is only one \erg\
component and $(\mN, \tnl, \mu_\mN)$ is \erg.  The \erg ity of both
$(\ta \times S^{d-1}, \phi^t, m)$ and $(\mD, T_\mD, \mu_\mD)$ follows
immediately.

As for the $K$-mixing property of the latter \dsy, once again we
construct its LSUMs as push-forwards or pull-backs of the LSUMs of,
say, $(\mN, \tnl, \mu_\mN)$. Here is where the hypothesis that $D$ is
made up of whole smooth boundary components comes into play: no
further cuts must occur during the push-forward/pull-back
process. Once a \hyp\ structure has been established for $(\mD, T_\mD,
\mu_\mD)$ one uses the general result of Pesin's theory \cite{p}
whereby the \sy\ decomposes into a countable number of positive-\me\
\erg\ components, over which a power of the map is piecewise
$K$-mixing \cite{ks}. Since in our case it is easy to prove that
$T_\mD^k$ is \erg\ for all $k \in \Z^+$ (the above-defined LSUMs are
LSUMs for any power of the map as well), this immediately implies that
the \sy\ is $K$-mixing.

\section{A few detailed arguments}
\label{sec-detail}

We have thus seen that, whenever Theorem \ref{thm-ltl} holds, the
desired result follows by fairly general arguments.  This final
section---which is the original core of this note---is devoted to
demonstrating that, under the stated assumptions, Theorem
\ref{thm-ltl} does hold for LTs in dimension 3 and higher.

The problem with \hyp\ \bi s in $d \ge 3$, as is common knowledge in
the field and perhaps not so common elsewhere, is that the so-called
\emph{fundamental theorem} (namely, local \erg ity; cf.\ Theorem
\ref{thm-ltl}\emph{(c)}) is not known for general semi-dispersing
tables, even in finite \me. The misinformation is due to the fact that
incorrect proofs of said theorem were believed valid until recently,
when B\'alint, Chernov, Sz\'asz and T\'oth \cite{bcst} pointed out the
mistake. In the same paper, these authors recover the proof for the
case of \emph{algebraic} Sinai \bi s, i.e., dispersing \bi s on the
torus with a finite number of scatterers, whose boundaries are made up
of a finite number or compact pieces of algebraic varieties.  To our
knowledge---if we exclude generic results that so far can claim no
specific examples \cite{bbt}---these are essentially the only
semi-dispersing tables, in $d \ge 3$, for which \erg ity is known.

Because the situation for Sinai \bi s is less than optimal, our
ability to prove local \erg ity for our models is also less than
optimal. In truth, we simply adapt the results of \cite{bcst} to the
framework at hand, much as a previous paper by one of us \cite{l1}
adapted the classical results on two-dimensional semi-dispersing \bi s
to planar Lorentz gases.

In what follows, we drop the subscript $\l$ from all the notation,
denoting our \dsy\ as $(\mM, T, \mu)$. To start with, one needs to
establish uniform \hyp ity, namely, statement \emph{(a)} of Theorem
\ref{thm-ltl}, because, as we shall see below, that is used in the
proof of \emph{(c)}.

In 2D, uniform \hyp ity descends rather easily from the fact that,
after a certain time past a dispersing collision (with curvature
bounded below), any \tr y has acquired a sufficient amount of \hyp
ity; namely, if one constructs an infinitesimal dispersing beam around
the \tr y, the beam has increased its dispersion by a large enough
factor, cf.~\cite{cls}.

For $d \ge 3$, a complication arises that is sometimes referred to as
the problem of \emph{astigmatism} \cite{b}.  It turns out that
dispersing beams of \tr ies that collide with a scatterer almost
tangentially (these are usually called \emph{grazing beams}) do not
acquire much additional dispersion.  This problem is circumvented if
one prescribes that, a positive percentage of the time, any given \tr
y undergoes a dispersing collision that has a non-negligible head-on
component \cite{bd}.
 
Assumption (A4) guarantees this. In fact, $\mM$ is defined as the
cross-section of all line elements on a dispersing boundary whose
outgoing velocity has a sufficiently large component along the normal
vector to the boundary. $\mM$ is thus a global cross section, with
corresponding map $T$. The same assumption also guarantees that, after
$K_3$ returns to $\mM$, a \tr y has traveled at least a distance $L$,
and thus (by (A3) as well) has acquired enough \hyp ity, relative to
the so-called orthogonal Jacobi metric in tangent space \cite{w,
  bd}. But, on $\mM$, because of the inequality in (\ref{def-malpha}),
this metric is equivalent to the natural Riemannian metric.  Therefore
$T^{K_3}$ is uniformly \hyp. But this implies the same for $T$, since
the orthogonal Jacobi metric is non-decreasing for tangent vectors
corresponding to dispersing beams.

In order to give a rigorous formulation of assumption (A6) and explain
how it is used in the proof of Theorem \ref{thm-ltl}\emph{(b)}, we
need to lay out some facts and a bit of extra notation.

It is common knowledge that \bi\ maps such as $T$ are discontinuous.
If $x \in \mM$ is the initial condition of a singular \tr y that has a
tangential collision or hits an edge before the next return to $\mM$,
then quite generally $x$ is a discontinuity point of $T$. We call such
$x$ a \emph{singular point} for the map $T$. (If $x$ is singular
because of a tangential collision, it can be seen that the
differential of $T$ blows up at $x$, whence the term `singular'.)

In our case, given the peculiar choice of $\mM$, cf.\
(\ref{def-malpha}), we must consider as singular also those \tr ies
which, after a collision, have a velocity $v$ such that
\begin{equation}
  \label{extra-si}
  v \cdot o_q = \eps.
\end{equation}
The corresponding line element $(q,v)$ is clearly a discontinuity
point for the Poincar\'e map of $\mM$.

Let $\si$ denote the set of all singular points of $T$ and define
$\sia := \si \cap \mMa$.  It is a well-known and easily derivable fact
that $\sia$ is decomposed into smooth portions of codimension-one
manifolds, each of which corresponds to a source of singularity (a
tangential scattering, an edge, or condition (\ref{extra-si}))
encountered before or at the next return to $\mM$; and to the
\emph{itinerary} of scatterers visited before that.  By (A2) and (A4),
the number of scatterers (and thus number of edges) that can be
visited before the next return to $\mM$ is bounded by a universal
constant. Therefore the number of smooth pieces that $\sia$ comprises
is also bounded by a universal constant.

Moreover, since the LT is algebraic in the sense of (A2), an easy
adaptation of the results of \cite{bcst} implies that $\sia$ is
actually a finite union of SSAVs, whose number is universally
bounded. (The proof of the algebraicity of the singularity set, in
\cite[\S5.1]{bcst}, does not use in an essential way that the
scatterer configuration is periodic there.) Notice that the extra
singularities due to (\ref{extra-si}) also give rise to
SSAVs. (Substitute eqn.\ (5.4) of \cite{bcst} with the polynomial
equation corresponding to (\ref{extra-si}).)

For $\delta > 0$, define
\begin{equation}
  \sia^{[\delta]} := \rset{x \in \mMa} {\mathrm{dist} 
  (x, \sia) < \delta}.
\end{equation}
The \me s of these neighborhoods play a pivotal role in the proof of
the \hyp\ properties of \bi s. The considerations in the previous
paragraph and the results of \cite[\S5.2]{bcst} imply that, as $\delta
\to 0$,
\begin{equation}
  \label{hd-10}
  \mathrm{Leb} ( \sia^{[\delta]} ) = \mathcal{O} (\delta),
\end{equation}
where $\mathrm{Leb}$ is the Lebesgue \me\ on $\mM$ (more precisely,
the Riemannian volume on $\mM$ corresponding to the distance
$\mathrm{dist}$; notice that $\mu$ is absolutely continuous w.r.t.\
$\mathrm{Leb}$.) The implicit constant in the r.h.s.\ of (\ref{hd-10})
depends in general on $\l$ and $\alpha$: we require the bound to be
uniform. More precisely we reformulate:
\begin{itemize}
\item[(A6)] There exists a universal constant $K_4 > 0$ such that, for
  all sufficiently small $\delta$, $\mathrm{Leb} ( \sia^{[\delta]} )
  \le K_4 \, \delta$.
\end{itemize}

By (\ref{hd-10}) it is not hard to generate examples of LTs satisfying
(A6). For example, one can start with a periodic algebraic LT and then
perform a (quenched \ra) \emph{algebraic perturbation}. By this we
mean that the equations of the perturbed scatterers are polynomials
whose coefficients are very close to the corresponding coefficients
for the unperturbed cell.  Another easy example is the case where
$\Omega$ is finite. In that case, say that $\bo_\alpha$ is contained
in the cell $C_n$. Then, by (A4), $\sia$ is completely determined by 
the local configurations in the cells $C_k$, with $n - K_5 \le k \le n +
K_5$ ($K_5$ being a universal constant). But there are only finitely
many possibilities for these local configurations, therefore (A6) is
implied by (\ref{hd-10}).

The argument that proves Theorem \ref{thm-ltl}\emph{(b)} is virtually
the same as in Lem.~3.2 of \cite{l1}. It is based on the old principle
that, in good \hyp\ \bi s, the exponential expansion of the LUMs, or
candidates therefor, is the dominant effect, compared to the cutting
operated by the singularities. Thus, locally along a given \o, one has
all the ingredients of Pesin's theory to prove the existence of a
LSM. The same ingredients then guarantee the absolute continuity of
the correspondong foliation.

More specifically, we want to show for $\mu$-a.a.\ $x \in \mM$,
a constant $C_0 = C_0(x)$ can be found such that
\begin{equation}
  \label{hd-20}
  \mathrm{dist} \! \left( T^{-k}x, \si \cup \partial \mM \right) 
  \ge C_0 \, k^{-3},
\end{equation}
for all positive integers $k$. Let us fix $x \in \mMa$. By the
finite-horizon condition, $T^{-k}x$ can only belong to a limited
portion of the phase space, namely $\bigcup_{\beta \in \mathcal{V}_k}
\mM_\beta$, where $\mathcal{V}_k = \mathcal{V}_k (\alpha)$ is the
index set of all the boundaries that can be visited within time $k$ by
a \tr y starting in $\bo_\alpha$. Using (A2) as well,
\begin{equation}
  \label{hd-25}
  \# \mathcal{V}_k \le K_6 \, k,
\end{equation}
for some $K_6 > 0$. So we have that the statement
\begin{equation}
  \label{hd-30}
  \mathrm{dist} ( T^{-k}x, \si \cup \partial \mM ) \le k^{-3}
\end{equation}
is equivalent to the statement
\begin{equation}
  \label{hd-40}
  x \in T^k \! \left( \bigcup_{\beta \in \mathcal{V}_k} \left( 
  \si_\beta \cup \partial \mM_\beta \right)^{[k^{-3}]} \right).
\end{equation}
By the invariance of $\mu$, (\ref{hd-25}) and (A6), the \me\ of the
r.h.s.\ of (\ref{hd-40}) is bounded by a constant times $k^{-2}$,
which is a summable series in $k$. By Borel-Cantelli applied to the
finite-\me\ space $(\mMa, \mu)$, the event (\ref{hd-40}), equivalently
(\ref{hd-30}), may happen infinitely often in $k$ only for a
negligible set of $x$, whence (\ref{hd-20}). The same reasoning
applies of course to every $\alpha$.

Finally, for the statement \emph{(c)}, we prove our version of the
fundamental theorem (cf.\ \cite[Sec.~4]{l1}) using the technique of
\emph{regular coverings}, as in \cite{kss} or \cite{lw}. This
technique requires a global argument (i.e., an estimate on objects
outside of the neighborhood under consideration) in one part only, the
so-called \emph{tail bound}.  The rest of the proof is local, thus
unable to distinguish between a finite- and an infinite-\me\ \bi: all
the standard arguments---including the exacting ones where one uses
that the singularity set and its images via the map are made up of
SSAVs \cite{bcst}---apply there.

The tail bound is the following statement: For all $x_0 \in \mM$, 
there exists a neighborhood $U$ of $x_0$ such that
\begin{equation}
  \label{hd-50}
  \mu \! \left( \rset{x \in U} { \mathrm{dist}_{W^s} \! \left (x, 
  \bigcup_{k>m} T^{-k} \si \right) < \delta} \right) = \delta \, o(1),
\end{equation}
as $m \to \infty$. Here $\mathrm{dist}_{W^s} (x, \cdot)$ is the
Riemannian distance along $W^s(x)$. (Compare (\ref{hd-50}) with the
statement of Lem.~4.4 of \cite{l1}, noticing that here we use $T^{-1}$
and $\si$, instead of $T$ and $\si^-$, the latter denoting the
singularity set of $T^{-1}$.)  Once again, there is no loss of
generality in choosing $U \subset \mMa$. By the earlier reasoning, the
only singularities whose images via $T^{-k}$ can get close to $x \in
\mMa$ are those in $\bigcup_{\beta \in \mathcal{V}_k}
\si_\beta$. Therefore (\ref{hd-50}) descends from the estimate:
\begin{eqnarray}
  && \mu \! \left( \rset{x \in \mMa \!} { \mathrm{dist}_{W^s} \! \left 
  (x, \bigcup_{k>m} T^{-k} \, \bigcup_{\beta \in \mathcal{V}_k} \si_\beta
  \right) < \delta} \right) \nonumber \\
  &=& \mu \! \left(  \bigcup_{k>m} \rset{x \in \mMa \!} { 
  \mathrm{dist}_{W^s} \! \left (x, \, T^{-k} \! \bigcup_{\beta \in 
  \mathcal{V}_k} \si_\beta \right) < \delta} \right) \nonumber \\
  &\le& \mu \! \left(  \bigcup_{k>m} \rset{x \in \mM \!} 
  { \mathrm{dist}_{W^s} \! \left ( T^k x, \bigcup_{\beta \in 
  \mathcal{V}_k} \si_\beta \right) < \delta c \lambda^k} \right) 
  \nonumber \\
  &\le&  \sum_{k=m+1}^\infty \mu \! \left( T^{-k} \! \left( 
  \bigcup_{\beta \in \mathcal{V}_k} \si_\beta^{[\delta c \lambda^k]} 
  \right) \right) \nonumber \\
  &\le& \mathrm{const} \ \delta \! \sum_{k=m+1}^\infty k \lambda^k. 
  \label{hd-60}
\end{eqnarray}
In the first inequality we have used the uniform \hyp ity of $T$
($\lambda<1$ is the contraction rate and $c$ is a suitable
constant). The third inequality follows from the invariance of $\mu$
and (A6).

\appendix

\section{Appendix: An example of a 3D Lorentz tube}
\label{app}

In this appendix we present an example of a three-dimensional LT that
verifies all the assumptions of Section \ref{sec-pre}.

We begin by constructing a template cell with its set of
scatterers. Then we check that the geometric assumptions (A2)-(A4)
hold for this cell. Finally, we specify how to use that to construct a
quenched \ra\ LT that verifies all assumptions.

The cell is a rectangular parallelepiped with a square base of side
length 1 and with height $h>1$. For the sake of the description, we
embed this solid in the $(x,y,z)$-space: the square base lies in the
$(y,z)$-plane, so that the edges of length $h$ are in the
$x$-direction. We call the $x$-direction `longitudinal' and any
direction orthogonal to it `transversal'.

Puncture each square facet with a square hole of side $g$, centered in
the center of the facet. These holes take the role of the gates $G^i$.
Populate the cell with the following scatterers: For each of the four
edges of length $h$, consider the cylinder of radius $\rho$ which has
that edge as its axis. Choose $1/2 < \rho < (1-g) / \sqrt{2}$, so that
the cylinders do not obscure the gates and intersect each other with a
positive angle (Fig.~\ref{ex_fig}). The portions of the cylinders that
are inside the cell are morally our scatterers. However, in order to
satisfy assumption (A3), we modify these quarter-cylinders by adding a
small positive curvature along the longitudinal direction. We do so in
a way that in every transversal section we still have four
quarter-circles whose radii satisfy the inequalities given
earlier---in other words, we have a so-called \emph{diamond \bi}; cf.\
Fig.~\ref{sect}. For lack of a better name, we call the resulting
solids `cigars'.

\figmpng{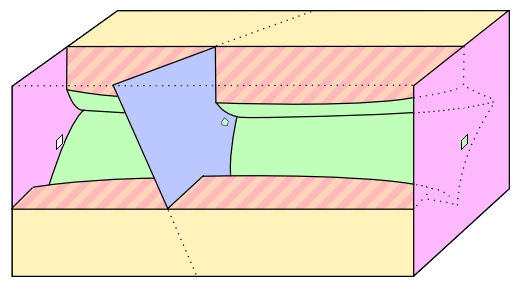} {9cm} {The template cell for our example of an LT
  satisfying (A1)-(A6).}

Moreover, we insert a `bulkhead' in the middle of the cell, that is,
we add an infinitesimally thin scatterer given by the intersection of
the plane passing through the center of the cell and orthogonal to the
vector $(1,1,1)$, and the cell itself. We punch an off-center, small,
polygonal hole through the bulkhead, so that it is possible to go from
one ``chamber'' to the other, but no free (i.e., collisionless) \tr
ies exist between the two gates of the cell.  Notice that by choosing
the parameter $h$ large enough, we can ensure that the maximum number
of consecutive collisions between the flat facets and the bulkhead is
3.

Next, we argue that this setting verifies (A4). We say that a \tr y
that hits a scatterer almost tangentially \emph{grazes} it. So
consider a \tr y that grazes a cigar. If this \tr y has a sizable
transversal component (compared to the longitudinal one) then its
projection on a transversal section is not far from a grazing \o\ in a
diamond \bi, which implies that the next bounce is not grazing. Let us
hence consider grazing \tr ies that are almost longitudinal. In this
case, within a certain (bounded above) time, the point will hit the
bulkhead in such a way that the next collision is necessarily
non-grazing and against a cigar, provided $h$ is not too small. On the
other hand, since the ``transversal angle'' at each intersection point
of two cigars is bounded below, no more that a bounded number of
collisions can be performed in said time. These arguments will be made
more quantitative below.

Finally, the full LT is made up of \ra\ algebraic perturbations (in
the sense of Section \ref{sec-detail}) of this local configuration. If
the perturbations are i.i.d.\ in each cell, (A1) is obviously
verified.  If they are sufficiently small, (A6) holds.  As for (A5),
it is not hard to see (by \erg ity of the inner dynamics in a cell, if
one will) that this assumption is verified as well.

\figmpng{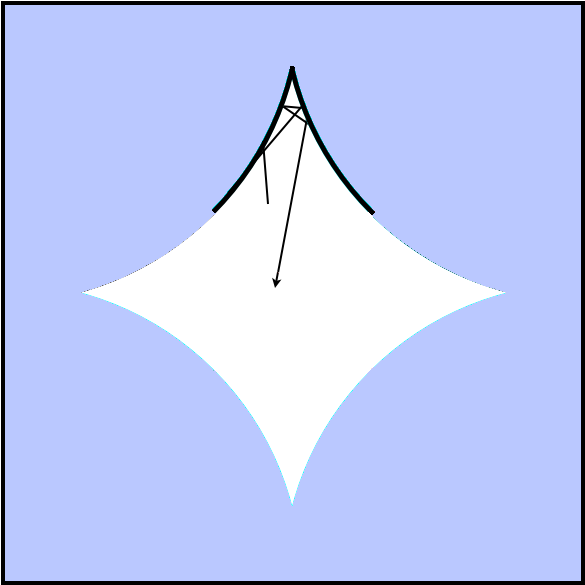} {6.5cm} {The longitudinal projection of a \tr y in the
  (purely cylindrical) cell, describing a planar \tr y in a diamond
  \bi. The indicated part of the boundary represents a `zone' (see
  text).}

At the request of an anonymous referee, we give estimates for the
universal constants $L$ and $K_3$ for this LT.  We begin by
considering a periodic tube made up of template cells that---contrary
to what we have imposed earlier---have no longitudinal curvature. So
the dispersing scatterers are truly cylindrical and, if we neglect the
bulkheads, describe the same diamond in each transversal section
(Fig.~\ref{sect}). Call $\gamma$ the angle at each vertex of the
diamond (so $\gamma = \arccos(1/2\rho)$). Fix $h > 10$, say, and $g <
10^{-1}$.

Denote by $v_x$ and $v_{yz}$, respectively, the longitudinal and
trasversal components of the velocity $v$ of the material point (thus,
$|v_x|^2 + |v_{yz}|^2 = 1$). The longitudinal projection of a piece of
\tr y that does not hit a bulkhead is a (planar) \tr y of the diamond
\bi. Consider a collision against a cylinder (for the 3D \bi) with
outgoing velocity $v$. Denoting by $\theta \in [0, \pi/2)$ the angle
of incidence of $v$ w.r.t.\ the normal at the collision point and by
$\varphi$ the corresponing angle for the projected \tr y, we have
$\cos \theta = |v_{yz}| \cos \varphi$.

Let us make a couple of observations on the dynamics in the diamond
\bi. To start with, we subdivide its boundary in 4 `zones', each zone
being defined as the set of all the boundary points that are closer to
a given vertex than to any other (see Fig.~\ref{sect}). Let $M :=
\lceil \pi/\gamma \rceil$, i.e., $M$ is the minimum integer $\ge
\pi/\gamma$. It can be seen that any \tr y in the diamond can have at
most $M$ consecutive collisions in the same zone, and at least one of
them will have an angle of incidence $\varphi < (\pi - \gamma)/2$ (to
see this, just ``linearize'' the corners). Also, the time between the
first collision in a zone and the first collision in another is
bounded below by $L_1 := 1/\sqrt{2} - \rho$ (this is the semilength of
the biggest square inscribed in the diamond).

Coming to three-dimensional \tr ies, we distinguish two cases: the \tr
ies with a good transversal component, defined by $|v_{xy}| \ge
10^{-2}$, and the other ones, which we call almost longitudinal.
Starting with either case, we want to find an upper bound for the time
before the next ``head-on'' collision---the parameter $\eta$ that
defines ``head-on'', cf.\ (A4), will de determined later. By
time-reversibility, twice this upper bound will be a good estimate for
$L$.

If a \tr y keeps a good transversal component during its entire visit
in a single zone, we know from above that there will be a collision
with angle of incidence $\theta \le \arccos( 10^{-2} \sin (\gamma/2) )
=: \pi/2 - \eta$, in a time less than $L_2 := 2h M$ ($2h$ is an
estimate for the horizon of our \bi). If not, we eventually fall in
the next case, which we consider right away.  An almost longitudinal
\tr y will remain such until it hits a bulkhead, and this occurs
necessarily within a time $L_3 := 3h (1 - 10^{-4})^{-1/2}$. The next
collision after that, given the inclination of the bulkhead and that
$h > 10$, is against a cylinder, with $\theta < \pi/2 - \eta$.

So, $L$ can chosen to be $2(L_2 + L_3)$.  During this time, the point
can visit at most $\lceil L/L_1 \rceil$ zones, therefore an upper
bound for the number of collisions it can have is $K_3 := \lceil L/L_1
\rceil M + \lceil 3L/h \rceil$ (the first term estimates the
collisions against the cylinders and the second term the collisions
against the flat boundaries).

\skippar

We treat the original LT as a perturbation of the one just considered.
If the radius of curvature of the cigars in the longitudinal direction
is, say, larger than $100h$, multiplying the above estimates by 100 
will certainly work.

\end{document}